\documentclass[12pt]{amsart}
\usepackage{amssymb}

\usepackage{graphicx}

\usepackage[]{amsfonts}

\def\underset#1#2{\mathrel{\mathop{\kern0pt #2}\limits_{#1}}}

\def\overset#1#2{\mathrel{\mathop{\kern0pt #2}\limits^{#1}}}

\def\couleur(#1 #2 #3)
	{
\special{ps::[end]}}

\def\sqr#1#2{{\vcenter{\vbox{\hrule height.#2pt
             \hbox{\vrule width.#2pt height#1pt \kern#1pt
             \vrule width.#2pt}
             \hrule height.#2pt}}}}

\def\st{\mathinner{\mkern1mu\raise1pt\hbox{.}				
		   \mkern1mu\raise4pt\hbox{.}
		   \mkern1mu\raise1pt\hbox{.}
		 }
         }

\def\bx#1{\setbox1=\hbox{\kern3pt{#1}\kern3pt}				
 \dimen1=\ht1 \advance\dimen1 by 3pt \dimen2=\dp1 \advance\dimen2 by 3pt
 \setbox1=\hbox{\vrule height\dimen1 depth\dimen2\box1\vrule}%
 \setbox1=\vbox{\hrule\box1\hrule}%
 \advance\dimen1 by .4pt \ht1=\dimen1
 \advance\dimen2 by .4pt \dp1=\dimen2 \box1\relax}

\def\k#1{\kern#1em}
\def\vci{\vrule  width.02em height1.47ex depth-.0ex}				
\def\11{{\rm\k{.2}\vci\k{-.37}1}}

\newtheorem{Theorem}{Theorem}[section]
\newtheorem{Proposition}[Theorem]{Proposition}

\begin{document}
\title{A mathematical model for Tsunami generation\\using a conservative 
velocity-pressure hyperbolic system}
\author{Alain-Yves LeRoux}
\address{UNIVERSITE BORDEAUX1, Institut Math{\'e}matiques de Bordeaux, 
UMR 5251,351,Cours de la Lib{\'e}ration, 33405, Talence Cedex}
\email{Alain-Yves.Leroux@math.u-bordeaux1.fr\\}
\maketitle
\begin{abstract} {
By using the Hugoniot curve in detonics as a Riemann invariant of a velocity-pressure 
model, we get a conservative hyperbolic system similar to the Euler equations. 
The only differences are the larger value of the adiabatic constant ($\gamma =8.678$ 
instead of $1.4$ for gas dynamics) and the mass density replaced by a 
strain density depending on the pressure. The model is not homogeneous 
since it involves a gravity and a friction term. After the seismic wave 
reaches up the bottom of the ocean, one gets a pressure wave propagating 
toward the surface, which is made of a frontal shock wave followed by 
a regular decreasing  profile. Since this regular profile propagates 
faster than the frontal shock waves, the amplitude of the pressure wave 
is strongly reduced when reaching the surface. Only in the case of a 
strong earth tremor the residual pressure wave is still sufficient to 
generate a water elevation with a sufficient wavelengths enable to propagate 
as a SaintVenant water wave and to become a tsunami when reaching the 
shore. We describe the construction of the model and the computation 
of the wave profile and discuss about the formation or not of a wave. 
  \ \par
{\hskip 1.5em}}\end{abstract}
We propose a model using a constant mass density because variable mass 
density models are often unstable, since a tiny variation of the density 
allways causes a large variation of the pressure. In the model presented 
here, the transported variable is not the mass density but a new variable, 
called the strain density, which has the same properties of conservation 
without this drawback of numerical unstability.\ \par
{\hskip 1.8em}The first section deals with the construction of the velocity-pressure 
model and the new conservative variable of strain density is designed 
in the second section. The new non homogeneous Euler-like model is then 
studied in the third section and one dimension numerical computations 
of the profile of the wave are reported in the fourth section, for different 
values of the friction coefficient. We conclude by some discussions about 
the emergence of a tsunami wave or not.\ \par
\section{The velocity-pressure model{\hskip 1.8em}}
\setcounter{equation}{0}The Hugoniot curves (see [6]) correspond to the 
linkage between the velocity $w$ and the pressure $p$ in a shock wave 
travelling through a material after an impact. By starting from a position 
at rest, the algebraic equation of such curves has the form
\begin{displaymath} 
p\ =\ \rho _{0}\ \left({\ c_{0}\ w\ +\ S_{0}\ w^{2}}\right) \ ,\end{displaymath} 
where $\rho _{0}$ is the mass density of the material, $c_{0}$ is its 
speed of sound, and $S_{0}$ is a dimensionless constant retated to this 
material, obtained from experiments. For water,  $\rho _{0}=1000\ kg\ m^{-3},\ c_{0}=1647\ m\ s^{-1}$ 
and $S_{0}=1.921.$ By solving with respect to $w$ , we get
\begin{displaymath} 
w(p)\ =\ \pm \ \alpha _{0}\ \left({\ \sqrt{1+\beta _{0}p}\ -\ 1\ }\right) 
\ ,\end{displaymath}  where $\alpha _{0}=\frac{c_{0}}{2\ S_{0}}=429\ m\ s^{-1}$ 
and $\beta _{0}=\frac{4\ S_{0}}{\rho _{0}\ c_{0}^{2}}=2.84\ 10^{-9}$ 
whose unit is the inverse of an energy. We notice the relation $\alpha _{0}\ \beta _{0}\ \rho _{0}\ c_{0}=2\ .$ 
The values of these parameters come from [5]. \ \par
{\hskip 1.8em}A general one dimension velocity-pressure model, with a 
constant mass density $\rho _{0}$, is made of a dynamics equation of 
the form 
\begin{equation} 
\rho _{0}\ \displaystyle \left({\displaystyle \ w_{t}+ww_{z}}\right) \ 
+\ p_{z}\ =\ 0\ \label{tsunami1}
\end{equation}  where $t$ is the time and $z$ is the position along a 
vertical upwards oriented axis (with $z=0$ at sea level), and a Hooke 
law of the form 
\begin{equation} 
p_{t}\ +\ w\ p_{z}\ +\ \rho _{0}\ c(p)^{2}\ w_{z}\ =\ 0\ \label{tsunami2}
\end{equation} where $c(p)>0\ $stands for the pressure depending wave 
velocity, to be identified. These two equations compose a hyperbolic 
system whose Riemann invariants have the form$\displaystyle \ w'(p)\ =\ \pm \ \frac{\displaystyle 1}{\displaystyle \rho _{0}\ c(p)}\ 
\ .$ By derivating the expression of $w(p)$ from the Hugoniot curves, 
one gets 
\begin{displaymath} 
w'(p)\ =\ \pm \ \frac{\alpha _{0}\ \beta _{0}}{2\ \sqrt{1+\beta _{0}p}}\ 
\ .\end{displaymath}  Identifying the two expressions and using $\alpha _{0}\beta _{0}\rho _{0}c_{0}=2$ 
lead to the formula
\begin{displaymath} 
c(p)\ =\ c_{0}\ \sqrt{1+\beta _{0}p}\ \ .\end{displaymath} This reads 
like a state law for our velocity pressure model. \ \par
{\hskip 1.8em}Up to now we were only concerned with homogeneous equations 
since Riemann invariants only exists in this case. We also have to take 
in account the gravity effects, since the pressure is increasing with 
the depth of water, and some friction effect since stillness is a stable 
configuration. We use a friction term of the Strickler type as usual 
in hydraulics, but other choices are possible and will lead to similar 
results. That way, the dynamics equation is replaced by
\begin{equation} 
\rho _{0}\ \displaystyle \left({\displaystyle \ w_{t}+ww_{z}}\right) \ 
+\ p_{z}\ +\ \rho _{0}\ g\ +\ k\ \displaystyle \left\vert{\displaystyle 
w}\right\vert w\ =\ 0\ ,\ \label{tsunami3}
\end{equation} where $g$ is the gravity constant and $k$ the friction 
parameter. The Hooke law and the state law are unchanged, and we take 
$g=9.8\ m\ s^{-1}$in the numerical experiments. The size of the friction 
parameter is a priori unknown and will be discussed later.\ \par
\section{The conservative strain density{\hskip 1.8em}}
\setcounter{equation}{0}We look for a quantity $q=q(p)$ satisfying the 
transport equation
\begin{equation} 
q_{t}\ +\ \displaystyle \left({\displaystyle qw}\right) _{z}\ =\ 0\ .\label{tsunami9}
\end{equation} Since $q$ depends on $p$ we get
\begin{displaymath} 
q'(p)\ \left({\ p_{t}\ +\ w\ p_{z}\ }\right) \ +\ q(p)\ w_{z}\ =\ 0\ ,\end{displaymath} 
to be compared with the Hooke law. We get 
\begin{displaymath} 
\frac{q'(p)}{q(p)}\ =\ \frac{1}{\rho _{0}\ c(p)^{2}}\ =\ \frac{1}{\rho 
_{0}\ c_{0}^{2}\ \left({1+\beta _{0}p}\right) }\ ,\end{displaymath} which 
is easily solved and gives
\begin{displaymath} 
q(p)\ =\ \ \left({\ 1+\ \beta _{0}\ p\ }\right) ^{\frac{1}{\beta _{0}\rho 
_{0}c_{0}^{2}}}\ =\ \left({\ 1\ +\ \beta _{0}\ p\ }\right) ^{\frac{\alpha 
_{0}}{2\ c_{0}}}\ \ ,\end{displaymath}  for the choice $q(0)=1\ .$ The 
wave velocity can be written as a function of the strain density $q$ 
as 
\begin{displaymath} 
c(q)=c_{0}\ q^{\frac{c_{0}}{\alpha _{0}}}\ .\end{displaymath}  We check 
now the conservation of the momentum $m=qw\ .$ We get
\begin{displaymath} 
m_{t}\ +\ \left({mw}\right) _{z}\ +\ \frac{q}{\rho _{0}}\ p_{z}\ +\ gq\ 
+\ \frac{k\ q}{\rho _{0}}\ \left\vert{w}\right\vert w\ =\ 0\ ,\end{displaymath} 
where 
\begin{displaymath} 
\frac{q}{\rho _{0}}\ p_{z}\ =\ \frac{2\ c_{0}}{\alpha _{0}\ \beta _{0}\ 
\rho _{0}}\ q^{2\frac{c_{0}}{\alpha _{0}}}\ q_{z}\ =\ c_{0}^{2}\ q^{2\frac{c_{0}}{\alpha 
_{0}}}\ q_{z}\ =\ c(q)^{2}\ q_{z}\ .\end{displaymath} We introduce a 
pressure term, standing as a strain pressure, 
\begin{displaymath} 
P(q)\ =\ c_{0}^{2}\ \frac{q^{2\frac{c_{0}}{\alpha _{0}}+1}}{2\frac{c_{0}}{\alpha 
_{0}}+1}\ =\ \frac{c_{0}^{2}}{\gamma _{0}}\ q^{\gamma _{0}},\ \ with\ 
\gamma _{0}=2\frac{c_{0}}{\alpha _{0}}+1\ (=\ 8.678)\ \end{displaymath} 
and we get the conservative equation for the momentum
\begin{equation} 
m_{t}\ +\ \displaystyle \left({\displaystyle \ mw\ +\ P(q)\ }\right) _{z}\ 
+\ gq\ +\ \frac{\displaystyle k}{\displaystyle \rho _{0}}q\displaystyle 
\left\vert{\displaystyle w}\right\vert w\ =\ 0\ .\label{tsunami8}
\end{equation}  The system made of ~(\ref{tsunami9}) ~(\ref{tsunami8}) 
has the form of the well known Euler equations for gas dynamics with 
a larger adiabatic coefficient $\gamma _{0}=8.678$ instead of the usual 
value $1.4$ for gases, and can be handled in the same way, especially 
for the shock waves. The same Rankine Hugoniot condition is valid, connecting 
the velocity of a shock wave to the two states $(q_{1},w_{1})$ and $(q_{2},w_{2})$ 
by
\begin{equation} 
z'(t)\ =\ \frac{\displaystyle w_{1}+w_{2}}{\displaystyle 2}\ +\ \frac{\displaystyle 
q_{1}+q_{2}}{\displaystyle 2}\ \displaystyle \sqrt{\displaystyle \frac{\displaystyle 
1}{\displaystyle q_{1}q_{2}}\ \frac{\displaystyle P(q_{2})-P(q_{1})}{\displaystyle 
q_{2}-q_{1}}}\ .\label{tsunami4}
\end{equation}  \ \par
\section{The profile of the strain wave{\hskip 1.8em}}
\setcounter{equation}{0}We first compute the state at rest. In case of 
stillness the equations ~(\ref{tsunami9}) and ~(\ref{tsunami8}) reduce 
to 
\begin{displaymath} 
q_{t}=0\ \ \ ,\ \ \ c_{0}^{2}\ q^{2\frac{c_{0}}{\alpha _{0}}-1}\ q_{z}\ 
+\ g\ =\ 0\ ,\end{displaymath} since $q{\not =}0.$ By denoting $q=q_{0}(z)$ 
the strain density at rest, the integration gives
\begin{displaymath} 
\frac{\alpha _{0}\ c_{0}}{2}\ q_{0}(z)^{2\frac{c_{0}}{\alpha _{0}}}\ +\ 
g\ z\ =\ Constant\ .\end{displaymath}  Recalling that $q^{2\frac{c_{0}}{\alpha _{0}}}=1+\beta _{0}p,$ 
we can use the atmospheric pressue $p_{a}$ at the surface ($z=0$) and 
get
\begin{displaymath} 
\frac{\alpha _{0}\ \beta _{0}\ c_{0}}{2}\ \left({\ p-p_{a}\ }\right) \ 
+\ g\ z\ =\ 0\ \ ,\ \ {\mathrm{that\ is\ }}\ \ \ p=p_{a}-\rho _{0}gz\ 
\ {\mathrm{or\ }}q_{0}(z)\ =\ \left({\ 1\ +\ \beta _{0}\ \left({p_{a}-\rho 
_{0}gz}\right) }\right) ^{\frac{\alpha _{0}}{2\ c_{0}}}\ \end{displaymath} 
which corresponds to the geostrophic equilibrium state.\ \par
{\hskip 1.8em}To compute the strain density profile we use the deviation 
variable $\eta =q-q_{0}$ and look for linkage of the form $m\ =\ A\ \eta \ -\ B\ ,$ 
as in any $q-m-$system with a source term (see [4], or annex below in 
Section 5), where $A$ and $B$ are constant. Since stillness is reached 
for $\eta =0,$ we have $B=0$. Besides, since $q_{0}$ does not depend 
on $t$, we have 
\begin{displaymath} 
\eta _{t}\ +\ A\ \eta _{z}\ =\ 0\ ,\end{displaymath} which means that 
$A$ corresponds to the wave velocity, that we name the reference velocity, 
corresponding to a reference state $(q_{ref},q_{ref}w_{ref})$ such that 
$A=w_{ref}+c(q_{ref})$. We compute, for $q>q_{0}$ which is always expected, 

\begin{displaymath} 
w=\frac{m}{q}\ =\ A\ \frac{q-q_{0}}{q}\ ,\ m_{t}=A\eta _{t}=-A^{2}\eta 
_{z}\ ,\ m_{z}=A\eta _{z}\ ,\ q_{z}=\eta _{z}-\frac{gq_{0}}{c(q_{0})^{2}}\ 
\end{displaymath} to be introduce into ~(\ref{tsunami8}) which becomes
\begin{displaymath} 
\left[{-A^{2}+2A^{2}\frac{q-q_{0}}{q}+c^{2}-A^{2}(\frac{q-q_{0}}{q})^{2}}\right] 
\eta _{z}\ +gq-\frac{gq_{0}}{c(q_{0})^{2}}\left({c^{2}-A^{2}(\frac{q-q_{0}}{q})^{2}}\right) 
\ +\ \frac{kq}{\rho _{0}}\frac{A^{2}}{q^{2}}\left({q-q_{0}}\right) ^{2}=0\ 
.\end{displaymath} This equation reduces to
\begin{displaymath} 
\left[{c^{2}-A^{2}\frac{q_{0}^{2}}{q^{2}}}\right] \ \eta _{z}\ +\ gq\ 
+\ \frac{gq_{0}}{c(q_{0})^{2}}\left[{A^{2}\left({\frac{q-q_{0}}{q}}\right) 
^{2}-c^{2}}\right] \ +\ \frac{k\ A^{2}}{\rho _{0}\ q}\left({q-q_{0}}\right) 
^{2}\ =\ 0\ .\end{displaymath} By multiplying by $q^{2}$ and using $c=c_{0}\ q^{\frac{c_{0}}{\alpha _{0}}}$ 
we get
\begin{equation} 
\displaystyle \left({\displaystyle q^{2}c^{2}-q_{0}^{2}A^{2}}\right) \ 
\eta _{z}\ +\ gq^{3}\displaystyle \left({\displaystyle 1-\displaystyle 
\left({\displaystyle \frac{\displaystyle q}{\displaystyle q_{0}}}\right) 
^{\displaystyle 2\frac{\displaystyle c_{0}}{\displaystyle \alpha _{0}}-1}}\right) 
\ +\ \frac{\displaystyle gA^{2}q_{0}}{\displaystyle c(q_{0})^{2}}\ +\ 
\frac{\displaystyle k\ q\ A^{2}}{\displaystyle \rho _{0}}\ \displaystyle 
\left({\displaystyle q-q_{0}}\right) ^{2}\ =\ 0\ .\label{tsunami5}
\end{equation} Since $q=\eta +q_{0},$ this is a differential equation 
which can be integrated by using standard numerical methods. An increasing 
profile is expected. Since the reference velocity $A\ $is far larger 
than $c(q)$ the coefficient $q^{2}c^{2}-q_{0}^{2}A^{2}$ of $\eta _{z}$ 
is always negative in practice. The two last terms are always positive, 
and the friction term is the predominant one. The term $gq^{3}\left({1-\left({\frac{q}{q_{0}}}\right) ^{2\frac{c_{0}}{\alpha _{0}}-1}}\right) 
$ is always negative and is always balanced by the friction term when 
the friction coefficent $k$ is not too small.\ \par
{\hskip 1.8em}The value of $A$ is determined by the strength of the sismic 
wave at the bottom of the ocean, whose depth is denoted $z_{f}.$ This 
corresponds to a reference state $(q_{ref},m_{ref}).$ We have 
\begin{displaymath} 
A\ =\ w_{ref}+c_{ref}\ \ ,\ \ m_{ref}\ =\ q_{ref}w_{ref}\ =\ A\ \left({q_{ref}-q_{0}(z_{f})}\right) 
\ \ ,\ c_{ref}\ =\ c_{0}\ q_{ref}^{\frac{c_{0}}{\alpha _{0}}}\ .\end{displaymath} 
 We get
\begin{displaymath} 
A\ =\ A\ \left({1-\frac{q_{0}(z_{f})}{q_{ref}}}\right) \ +\ c_{ref}\ ,\end{displaymath} 
whichs gives $A\ =\ c_{ref}\ \frac{q_{ref}}{q_{0}(z_{f})}\ \ {\mathrm{and}}\ w_{ref}\ 
=\ c_{ref}\ \left({\frac{q_{ref}}{q_{0}(z_{f})}-1}\right) \ .$\ \par
{\hskip 1.8em}Now we can compute the profile of the stain wave as the 
solution of ~(\ref{tsunami5}).\ \par
Figure 1 presents a series of numerical computation tests using the reference 
value $q_{ref}=1.1296$, for a depth $z_{f}\ =\ 3700$ meters, that is 
an increasing of about 1500\% above the natural pressure on the bottom 
of the ocean. \ \par

\begin{figure}[h]
\begin{center}
\rotatebox{0}{\resizebox{10cm}{!}{\includegraphics{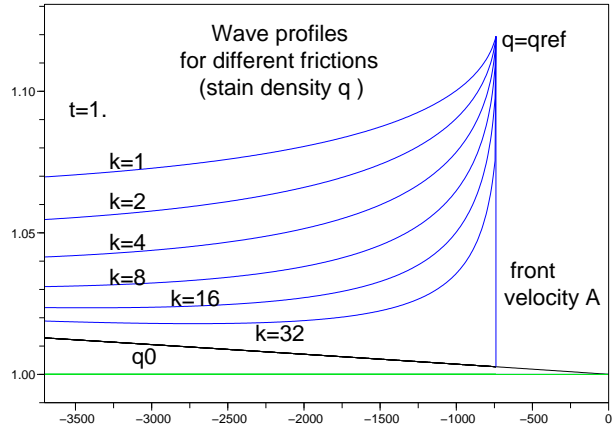}}}
\caption{The profiles depend on the friction coefficient}
\label{ts0}
\end{center}
\end{figure}
\ \par
{\hskip 2.1em}The velocity field  $w_{ref}$ increases from a few meters 
per second near the bottom to more than 300 meters per second near the 
front shock wave drawn here, which is here an hypothetic one. The tests 
performed with too small friction coefficients (observed here for $k<0.15$) 
yield decreasing profiles, as expected from the remark above about the 
size of the friction term. The real front shock wave will progress more 
slowly than the strain wave, and its amplitude will decrease rapidly, 
and the effective values of the velocity field $w_{ref}$ will be strongly 
reduced near the front schock wave.  \ \par
{\hskip 1.8em}The front shock wave connects the geostrophic equilibrium 
state $q_{0}$, with the velocity $w=0$, to a value $q$ on the strain 
wave, with the velocity $w=A\left({1-\frac{q_{0}}{q}}\right) .$ The Rankine 
Hugoniot condition ~(\ref{tsunami4}) gives the velocity of this shock 
wave, which reads here:
\begin{equation} 
z'(t)\ =\ \frac{\displaystyle A}{\displaystyle 2}\displaystyle \left({\displaystyle 
1-\frac{\displaystyle q_{0}}{\displaystyle q}}\right) \ +\ \frac{\displaystyle 
q_{0}+q}{\displaystyle 2}\ \displaystyle \sqrt{\displaystyle \frac{\displaystyle 
1}{\displaystyle q\ q_{0}}\ \frac{\displaystyle P(q)-P(q_{0})}{\displaystyle 
q-q_{0}}}\ \ ,\ \ {\mathrm{with}}\ P(q)\ =\ \frac{\displaystyle c_{0}^{2}}{\displaystyle 
\gamma _{0}}\ q^{\displaystyle \gamma _{0}}\ .\label{tsunami6}
\end{equation}  We have the following result:\ \par
\begin{Proposition}  The shock wave propagates slower than the stain 
wave, and faster than the local wave speed $c(q_{0}(z))$,  that is 
\begin{displaymath} 
c_{0}\ q_{0}(z(t))^{\frac{c_{0}}{\alpha _{0}}}\ <\ z'(t)\ <\ A\ .\end{displaymath} 
\ \par
\end{Proposition}
\textbf{Proof: } We fix $t$ and set $q_{s}(t)$ and $q_{0}=q_{0}(z(t))$, 
as the left and right values of the shock wave. Then the shock velocity 
reads
\begin{displaymath} 
z'(t)\ =\ \frac{A}{2}\ \frac{q_{s}-q_{0}}{q_{s}}\ +\ \frac{q_{s}+q_{0}}{2}\ 
\sqrt{\frac{1}{q_{s}q_{0}}\ \frac{P(q_{s})-P(q_{0})}{q_{s}-q_{0}}}\ .\end{displaymath} 
We shall use some $\xi \in ]q_{0},q_{s}[$ such that 
\begin{displaymath} 
\frac{P(q_{s})-P(q_{0})}{q_{s}-q_{0}}\ =\ P'(\xi )\ =\ c_{0}^{2}\ \xi 
^{\frac{2c_{0}}{\alpha _{0}}}\ .\end{displaymath}  We have
\begin{displaymath} 
z'(t)\ <\ A\ \ {\Longleftrightarrow}\ \frac{q_{s}+q_{0}}{2}\ \sqrt{\frac{1}{q_{s}q_{0}}\ 
\frac{P(q_{s})-P(q_{0})}{q_{s}-q_{0}}}\ <\ A\ \left({1-\frac{1}{2}\ +\frac{q_{0}}{2q_{s}}}\right) 
\ ,\end{displaymath}  where $\frac{q_{s}+q_{0}}{2}\ =\ q_{s}\left({\frac{1}{2}\ +\ \frac{q_{0}}{2q_{s}}}\right) 
\ .$ We get
\begin{displaymath} 
z'(t)\ <\ A\ \ {\Longleftrightarrow}\ q_{s}\left({\frac{1}{2}+\frac{q_{0}}{2q_{s}}}\right) 
\ \sqrt{\frac{1}{q_{s}q_{0}}\ c_{0}^{2}\xi ^{\frac{2c_{0}}{\alpha _{0}}}}\ 
<\ A\ \left({\frac{1}{2}\ +\frac{q_{0}}{2q_{s}}}\right) \ ,\end{displaymath} 
 which reduces to
\begin{displaymath} 
c_{0}\sqrt{\frac{q_{s}}{q_{0}}\ }\ \xi ^{\frac{c_{0}}{\alpha _{0}}}\ <\ 
A\ =\ c_{0}\ q_{ref}^{\frac{c_{0}}{\alpha _{0}}}\ \frac{q_{ref}}{q_{0}(z_{f})}\ 
,\end{displaymath}  and is equivalent to
\begin{displaymath} 
1\ <\ \left({\frac{q_{ref}}{\xi }}\right) ^{\frac{c_{0}}{\alpha _{0}}}\ 
\frac{q_{ref}}{q_{0}(z_{f})}\ \sqrt{\frac{q_{s}}{q_{0}}}\ ,\end{displaymath} 
which is true since $q_{ref}>q_{0}(z_{f})$ and $q_{ref}>q_{s}>\xi >q_{0}\ .$\ 
\par
{\hskip 1.8em}On the other hand, since the expression
\begin{displaymath} 
\ \frac{A}{2}\ \frac{q_{s}-q_{0}}{q_{s}}\ +\ \frac{q_{s}+q_{0}}{2}\ \sqrt{\frac{1}{q_{s}q_{0}}\ 
\frac{P(q_{s})-P(q_{0})}{q_{s}-q_{0}}}\ \end{displaymath}  is an increasing 
function of $q_{s}$ for $q_{s}\geq q_{0}$, we get obviously the other 
inequality.(End of proof)\ \par
\ \par
{\hskip 1.8em}Now we can construct the whole wave, made of a regular 
part corresponding to a part of the strain wave and a front shock whose 
position is determined by the Rankine Hugoniot condition ~(\ref{tsunami6}) 
intertpreted as a differential equation whose solution $z(t)$ gives the 
position of the shock. Figure 2 shows the different positions of the 
shock wave for different values of the friction coefficient. \ \par

\begin{figure}[h]
\begin{center}
\rotatebox{0}{\resizebox{10cm}{!}{\includegraphics{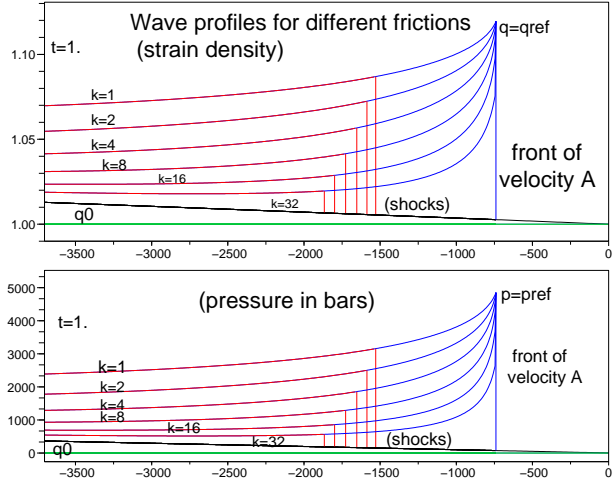}}}
\caption{The shock reduces the wave amplitude}
\label{ts1}
\end{center}
\end{figure}
\ \par
{\hskip 1.8em}The velocity $A$ is equal to $2932.5\ m/s$ that is $A=1.78\ c_{0}$. 
The shock on the bottom of the ocean, at the depth $z_{f}=-3700m,$ corresponds 
to 15 times the value of the usual geostrophic pressure, thai is  
\begin{displaymath} 
p_{bottom}\ =\ 363.6\ 10^{5}\ pascal\ =\ 363.6\ bars\ \ \ ,\ \ p_{ref}\ 
=\ 5.454\ 10^{8}\ pascal\ =\ 5.454\ kbars\ .\end{displaymath} This corresponds 
to strain densities 
\begin{displaymath} 
q_{fond}\ =\ 1.01288\ \ \ \ ,\ \ \ \ q_{ref}\ =\ 1.1296\ \ .\end{displaymath} 
 We notice that the variation from $q_{bottom}=q_{0}(z_{f})$ to $q_{ref}$ 
corresponds to an increasing of 11.5 \% which is an increasing of about 
1400 \% of the pressure. The velocity of the wave is computed from the 
values  
\begin{displaymath} 
c_{ref}\ =\ c_{0}\ q_{ref}^{\frac{c_{0}}{\alpha _{0}}}\ =\ 2629.5\ m/s\ 
\ \ ,\ w_{ref}\ =\ c_{ref}\ \frac{q_{ref}}{q_{f}}\ =\ 303\ m/s\ \ .\ 
\end{displaymath} We notice that the velocities $c_{ref}$ and $A=w_{ref}+c_{ref}\ $are 
for more important than $c_{0}$. The real profiles are drawn in red, 
and are shaped as a part of the strain wave profile cut by the slower 
shock wave. The eliminated parts are drawn in blue. We observe an important 
difference between the velocites of the different fronts, for several 
values of the friction coefficient.\ \par
\section{A wave or not ?{\hskip 1.8em}}
\setcounter{equation}{0}We denote by $H\ (=-z_{f}{\mathrm{\ here}})$ 
the mean depth of the ocean. The wavelength $\lambda $ of the Saint-Venant 
waves must satisfy a condition of the form
\begin{equation} 
\lambda \ \geq \ 2\ N\ \frac{\displaystyle H\ \displaystyle \sqrt{\displaystyle 
gH}}{\displaystyle c_{s}}\ \ ,\label{tsunami7}
\end{equation}  where $N$ the number of sonic interactions (back-and-forth) 
between the bottom and the surface of the ocean. This condition means 
that along a horizontal distance of a wavelength $\lambda $, there are 
at least $N$ such sonic interactions . The use of the Saint-Venant model 
is as more appropriate as $N$ is great. As in [2] for Rogue waves, we 
propose to require $N\geq 25,$ which implies for example a wavelength 
greater than  $21400\ m$ for an ocean depth of $3700m$. A tsunami wave 
is expected when the condition ~(\ref{tsunami7}) is fullfilled. A linear 
model for the surface elevation propagation was proposed in the historical 
paper [1] by K.Kajiura.\ \par

\begin{figure}[h]
\begin{center}
\rotatebox{0}{\resizebox{10cm}{!}{\includegraphics{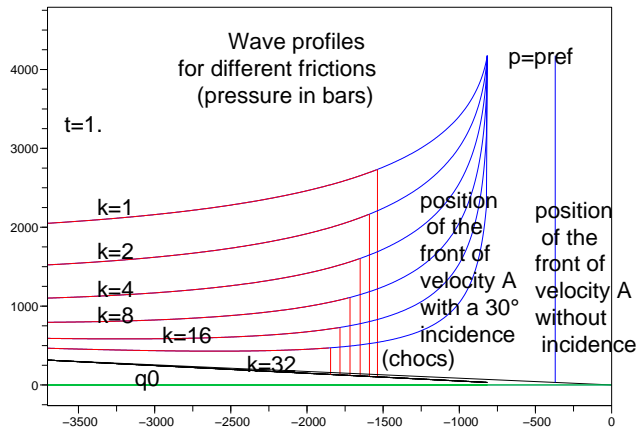}}}
\caption{The effect of the incidence}
\label{ts2}
\end{center}
\end{figure}
\ \par
{\hskip 1.8em}The numerical tests in section 3 show the effect of the 
friction, whose effect is to shape up the amplitude of the wave near 
the front shock. At the same time, this front shock magnitude is eroded 
progressively as the wave propagates upwards to the surface. Since the 
amplitude of the front shock decreases, the velocity of this front shock 
decreases too. By following the propagation of a wave with an incidence 
of angle $\phi $ (that is only changing $g$ into $g\ cos\phi $ and $z$ 
into $z/cos\phi $),  we get a longer path to travel with a weaker gravity 
constant) and a delayed wave compare to the case withou incidence. This 
is show on Figure 3.    \ \par
{\hskip 1.8em}The value of the initial amplitude has a capital effect. 
It must be large enough to get, after erosion by the friction, a remaining 
wave near the surface which is sufficient to raise up the sea surface 
and provoke a wave. The physical wave starts as a sperical wave, and 
propagates according to the incidences. The part with a small incidence 
will reach the surface later and will help the formation of the water 
wave. The part with a larger incidence will disappear because of the 
friction effect. The question of the value of the friction coefficient 
stays open, since the use of the strain density was never done before. 
It seems from the numerical tests that the correct values lay between 
1 and 10. Too large values provoke a sharp front wave which erodes rapidly 
and will never reach the surface with a sufficient amplitude to make 
a wave, which is not expected, since sometimes, tsunamis really occur. 
\ \par
\section{Annex: the source wave linearity{\hskip 1.8em}}
\setcounter{equation}{0}We consider a general 2x2 hyperbolic system whose 
first equation has the form
\begin{equation} 
q_{t}+m_{x}=0\ \ .\label{tsunami21}
\end{equation}  We denote by $\lambda _{1}$ and $\lambda _{2}\geq \lambda _{1}$ 
the eigenvalues of the flux matrix, which depend on $q$ and $m$ only. 
Then the general form of the second equation is 
\begin{displaymath} 
m_{t}\ +\ \left({\lambda _{1}+\lambda _{2}}\right) \ m_{x}\ -\ \lambda 
_{1}\lambda _{2}\ q_{x}\ =\ S(q,m)\ ,\end{displaymath}  or  
\begin{equation} 
m_{t}\ +\ 2\ u\ m_{x}\ +\ \displaystyle \left({\displaystyle c^{2}-u^{2}}\right) 
\ q_{x}\ =\ S(q,m)\ ,\label{tsunami22}
\end{equation}  by using the notations 
\begin{displaymath} 
u=\frac{\lambda _{1}+\lambda _{2}}{2}\ \ ,\ \ c=\frac{\lambda _{2}-\lambda 
_{1}}{2}\ \ .\end{displaymath}  and $S(q,m)$ is a source term, assumed 
to be not identically zero.  We have the following result:\ \par
\begin{Theorem} The nonlinear non homogeneous system ~(\ref{tsunami21}),~(\ref{tsunami22}) 
admits non constant local solutions which are also solutions to the linear 
homogeneous system 
\begin{equation} 
q_{t}\ +\ A\ q_{x}\ =\ 0\ \ ,\ \ m_{t}\ +\ A\ m_{x}\ =\ 0\ ,\label{tsunami23}
\end{equation} where $A$ is a real constant, with the linkage $Aq-m=B,$ 
another real constant.\ \par
\end{Theorem}
\ \par
\textbf{Proof: }We look for local solution with a linkage of the form 
$m=m(q)$. Then the system becomes
\begin{displaymath} 
q_{t}+m'(q)q_{x}=0\ \ ,\ \ m'(q)\ \left({q_{t}+2uq_{x}}\right) \ +\ \left({c^{2}-u^{2}}\right) 
q_{x}\ =\ S(q,m(q))\ .\end{displaymath} Since $q_{t}=-m'(q)\ q_{x}\ ,\ $and 
$S(q,m(q))\not\equiv 0$, the second equation becomes
\begin{equation} 
\frac{\displaystyle \displaystyle \left({\displaystyle u(q,m(q))-m'(q)}\right) 
^{2}-c(q,m(q))^{2}}{\displaystyle S(q,m(q))}\ q_{x}\ =\ 1\ ,\label{tsunami24}
\end{equation} which has the form 
\begin{displaymath} 
\psi '(q)\ q_{x}\ =\ 1\ ,\end{displaymath} by introducing a real function 
$\psi (q)$ whose derivative is$\displaystyle \ \psi '(q)\ =\ \frac{\displaystyle \displaystyle \left({\displaystyle 
u(q,m(q))-m'(q)}\right) ^{2}-c(q,m(q))^{2}}{\displaystyle S(q,m(q))}\ 
.$ Now, integrating with respect to $x\ $gives $\psi (q)=x-K(t)\ ,$ where 
$K(t)\ $is an integrating constant which may depend on $t$. Next, derivating 
with respect to $t$ gives $\psi '(q)\ q_{t}=-\ K'(t)\ ,$ where we set 
$q_{t}=-m'(q)\ q_{x}\ .$Hence we get $m'(q)\ \psi '(q)\ q_{x}=K'(t)\ ,$ 
and recalling that $\psi '(q)\ q_{x}=1\ ,$it remains $m'(q)=K'(t).$ A 
new derivation with respect to $x$ leads to $m''(q)\ q_{x}=0\ ,$ and 
since the solution is not constant, we get $m''(q)=0,\ $that is $m(q)=Aq-B,\ $with 
some constants $A$ and $B.$ Next, $K'(t)=m'(q)=A,$ and we have got $\displaystyle 
\psi (q)=x-At-x_{0}\ ,$ for some constant $x_{0},$ or, locally, $\displaystyle 
q=\psi ^{\displaystyle -1}(x-x_{0}-At)$, which satisfies to ~(\ref{tsunami23}), 
since $m_{t}+Am_{x}=m'(q)\ \left({q_{t}+Aq_{x}}\right) =0\ $.(end of 
proof)\ \par
\ \par
{\hskip 1.5em}This result is a very general one, since no special hypotheses 
were needed on the second equation ~(\ref{tsunami22}). Such waves are 
very common in the nature: water waves such as roll waves, rogue waves, 
tidal bore waves or also many other waves as reported in [2] or [4]. 
For example the double property of being either a solution to a non linear, 
non homogeneous systems and a linear homogeneous system provides the 
linkage between acoustics and gas dynamics in a wind instrument (see 
[3]). \ \par
{\hskip 1.8em}In the case of a conservative system, invariant by Galilean 
transform, the only choice of the function $u$ is reduced to $\displaystyle 
u(q,m)=m/q\ .$ We easily construct this way the usual Saint Venant system 
in hydraulics or the Euler equations in gas dynamics. \ \par
{\hskip 1.8em}In Section 3, the state at rest is not $q=0,\ m=0$, but 
$q=q_{0},\ m=0,$ so we look for a solution with the linkage $A\eta -m=B$ 
since $\eta =q-q_{0}=0$ at rest. We find a differentiel equation which 
is more complex than ~(\ref{tsunami24}). \ \par
\section{Biliography{\hskip 1.8em}}
\setcounter{equation}{0}We have used the idea of source term linearization 
effects on waves as in [2],[3] and [4]. Other models are developped in 
[1], [5] and [7]. The data [6] were used to valuate the parameters.\ 
\par
\ \par
[1] K.Kajiura, Tsunami source, energy and the directivity of wave propagation, 
Bulletin of the earthquake research institute, Vol.48, 835-869 (1970).\ 
\par
[2]A.-Y.LeRoux, M.-N.LeRoux, A Mathematical Model for Rogue Waves, using 
Saint Venant Equations with Friction. Conservation laws pr. server: www.math.ntnu.no$/$conservation$/$2005$/$048.html\ 
\par
[3] A.-Y.LeRoux, Sound traveling-waves in wind instruments as solutions 
to non linear non homogeneous gas dynamics equations.  www.math.ntnu.no$/$conservation$/$2008$/$010.html\ 
\par
[4]A-Y.LeRoux, M.-N.LeRoux,  Source waves-  www.math.ntnu.no$/$conservation$/$2004$/$045.html 
\ \par
[5] D.Dutykh, F.Dias, Y.Kervella. Linear theory of wave generation by 
a moving bottom , C. R. Acad. Sci. Paris, Ser. I 343, 499-504, (2006)\ 
\par
[6] Selected Hugoniot,  Report LA-4167-MS, Group GMX-6, Los Alamos Scientific 
Laboratory, University of California, Los Alamos, NM-87544 (1969).\ \par
[7] Y.Kervella, D.Dutykh, F.Dias. Comparison between three-dimensional 
linear and nonlinear tsunami generation models. Theoretical and Computational 
Fluid Dynamics, 21, 245-269, (2007).\ \par
\ \par
\ \par

\end{document}